\documentclass{article}
\usepackage[margin=1in]{geometry}

\usepackage[english]{babel}

\usepackage{graphicx}
\usepackage{hyperref}
\usepackage{color}

\usepackage{amssymb}
\usepackage{amsmath}
\usepackage{pifont}

\usepackage{float}
\usepackage{caption}
\usepackage{subcaption}
\graphicspath{{figures/}}

\usepackage{algorithm2e}
\DeclareMathOperator{\im}{im}

\title{A randomised non-descent method for global optimisation}
\author{Dmitry A. Pasechnyuk (MIPT, MBZUAI), Alexander Gornov (ISDCT SB RAS)}
\date{March 2023}

\begin{document}

\maketitle

\begin{center}
\begin{minipage}{0.65\textwidth}
    \begin{center}
    \textbf{Abstract}
    \end{center}
    \quad This paper proposes novel algorithm for non-convex multimodal constrained optimisation problems. It is based on sequential solving restrictions of problem to sections of feasible set by random subspaces (in general, manifolds) of low dimensionality. This approach varies in a way to draw subspaces, dimensionality of subspaces, and method to solve restricted problems. We provide empirical study of algorithm on convex, unimodal and multimodal optimisation problems and compare it with efficient algorithms intended for each class of problems.\newline
    
    
\end{minipage}
\end{center}

\section{Introduction}

Speaking of optimisation as a computational discipline, practitioner and therefore a theorist face with a great collection of functions with no structure sufficient to guarantee that their optimum can be found in reasonable time. The only properties of these functions which a designer of algorithms can rely on are continuity and boundedness of domain. This is what we mean by global optimisation. Most of natural problems in a non-simplified form satisfy this description, textbook example is intermolecular potential optimisation. We hope that mentioned properties are natural enough to cover most of real-life variational phenomena. Thus, solving global optimisation problems is essential purpose of applied mathematics.

Many practical general-purpose global optimisation methods represent a heuristic, which has proved the efficiency in applications. Great part of global optimisation as a science consist of systematic collection of such heuristics and computational facts about them. Classic works reflecting the described methodology include \cite{torn1989global,jones2001taxonomy}. To systematise numerical tests, the benchmarks are developed, for example \cite{gavana2018global}. Despite the growth of another, theoretical part of global optimisation science, this paper belongs to former one, and describes and tests a heuristic, which first appeared in preprint \cite{pasechnyuksolar}.

Global optimisation makes severe demands on a general-purpose heuristic. Continuous functions challenge us with expensive-to-evaluate function's value and gradient, curse of dimensionality, numerous spurious optima (local extrema and saddle points). Further moreover, descent direction given by anti-gradient ceases to carry global information about function if no growth conditions introduced. This is the reason why local descent-type methods are not widely used in practice of global optimisation as standalone algorithms, but usually serve as an auxiliary local minimisation procedure for basin-hopping-type algorithms \cite{wales1997global}.

The philosophy behind the proposed approach is following. Since the monograph \cite{nemirovskij1983problem}, the framework of local oracles is the dominant point of view on optimisation algorithms. In broad terms, in the beginning of algorithm's operation it has no information about the particular function to optimise and cannot classify it to point the optimum. Every next iteration of the algorithm calls an oracle, which gives some information to classify function more precisely. For example, if class of functions consist of two functions with no common values, one evaluation of function will fully characterise the function and its optimum will be known from pre-knowledge. When enough information is received to narrow the class of possible function so that their optima are close, optimisation completes. Further development of this concept for the case of local oracles is a foundation of convex optimisation. We cannot rely on local information, so we receive the information about the function through \textit{hypotheses}. Every next iteration of our algorithm test the hypothesis of some form. Hypothesis may be unsuccessful, or \textit{unproductive}, and give no information for optimisation (in this case, algorithm go to next step), or \textit{productive} and give some information to construct a step. When algorithm cannot come to a productive hypothesis for a long time, it indicates that good enough approximation of the optimum is already found. The particular instance, which is called Solar method, of the meta-algorithmic scheme described above is proposed in this paper. Its hypotheses are of the form: ``Does the point with optimal function's value in section of feasible set by a given subspace improves current best function's value?'' --- and its step is a jump to that point if hypothesis is verified.

The paper is organised as follows. In Section~\ref{sec:review} we summarise the development of algorithms allied to one described in this paper and provide a literature review. The algorithm itself is described in Section~\ref{sec:algorithm}, in Section~\ref{sec:subspace} this description is extended with description of ways to draw a random subspace. Section~\ref{sec:numerical} consists of empirical evaluation of algorithm's efficiency on test functions representing convex, unimodal and global optimisation problems in Sections~\ref{sec:convex}, \ref{sec:unimodal} and \ref{sec:global} correspondingly. There we compare proposed algorithm with algorithms known to be efficient on corresponding class, and test the dependencies of algorithm's convergence on parameters of problem and algorithm's hyperparameters. 

\section{Related work} \label{sec:review}

One forerunner of the algorithm we propose is coordinate descent method. It is known to be especially efficient in large-scale optimisation problems. In convex optimisation setting it provably converges to the optimum with known convergence rates. Classic works devoted to the analysis of coordinate descent methods are \cite{nesterov2012efficiency,richtarik2014iteration}. At the same time, this method is not widely used in non-convex and global optimisation. Attempts to go beyond the class of convex problems were made \cite{patrascu2015efficient,xu2017globally}, but corresponding results are as well of a local nature.

Another branch related to optimisation over random subspaces is sketching. This technique consists in implicit choice of a subspace by proper preconditioning of the gradient. Unlike coordinate descent method, this approach is more general in variants of subspace to choose. It became popular in community developing randomised optimisation algorithms \cite{gower2021stochastic,grishchenko2021proximal,kozak2021stochastic} and has acquired extensions for second-order oracle \cite{gower2019rsn} and non-convex problems \cite{cartis2022randomised}. Unlike the algorithm proposed in this paper, efficacy of these methods is based on Johnson--Lindenstrauss lemma or similar theoretical principles.

Second closest predecessor of our algorithm is classical steepest gradient descent with line-search. Properties of related methods in global optimisation environment is reflected in \cite{zhigljavsky2012theory} to some extent. However, algorithms with global line-search along descent or random directions are out of mainstream and it is hard to find modern works on this topic, as is easily seen if trying to browse in Google Scholar.

\section{The algorithm} \label{sec:algorithm}

One iteration of the algorithm, which is called Solar method, consists of three principal steps: construction of set to restrict the problem (random drawn subspace or manifold), solving restricted optimisation problem, and jump to the new point if it improves current record.

In general case, random manifold can be constructed using cloud of $< p$ reference points to interpolate or approximate through them. In dependence on method, one should store $p$ points in some data structure with cheap insertion and extraction of minimum of finding closest neighbours. In variants we consider it is necessary to maintain $p$ minimums, which can be done by heap or simple sorted list. If closest neighbours need to be found for construction, one should use $k$-d tree.

Pseudocode of Solar method is listed in Algorithm~\ref{alg:solar}. Naming of the randomly generated inclusion $r$ explains the name of method: the metaphor is that we move along random rays from current sun-like point.

\RestyleAlgo{ruled}
\begin{algorithm}
    \caption{Solar method}
    \label{alg:solar}
    
    \KwData{Number of outer iterations $K \in \mathbb{N}$, number of total inner iterations $N \in \mathbb{N}$, total dimensionality $n \in \mathbb{N}$, number of base variables $b \in \mathbb{N}: 1 \leq b < n$, number of probes $p \in \mathbb{N}$, convex indicator function $\chi: \mathbb{R}^n \to \{0, +\infty\}$ of feasible set $Q \in \mathbb{R}^n$, initial point $x \in Q$, zeroth-order oracle $f: \mathbb{R}^n \to \mathbb{R}$[, first-order oracle $\nabla f: \mathbb{R}^n \to \mathbb{R}^n$]}
    \KwResult{$x_\text{best} = T.\text{extract\_min}() \in \mathbb{R}^n$}

    Initialise data structure $T$\;
    $T.\text{insert}((f(x), x))$\;
    
    \For{$i = 1, ..., K$}
    {
        Choose $|B| = b$ random unique indices from set $\{1, ..., n\}$\;
        \For{$j = 1, ..., \lfloor N / K \rfloor$}
        {
            \For{$k = 1, ..., p$}
            {
                $(f_k, x_k) = T.\text{extract\_min}()$\;
            }
            Construct the ray $r: \mathbb{R}^b \hookrightarrow \mathbb{R}^n$ using points $(x_1, ..., x_p)$[ and gradients $(\nabla f(x_1), ..., \nabla f(x_p))$]\;
            $\displaystyle x_\text{candidate} \in \arg \min_{x \in \im r} \{f(x) + \chi(x)\}$\;
            $T.\text{insert}((f(x_\text{candidate}), x_\text{candidate}))$\;
            \For{$k = 1, ..., p-1$}
            {
                $T.\text{insert}((f_k, x_k))$\;
            }
        }
    }
\end{algorithm}

\subsection{Subspace, or manifold, sampling} \label{sec:subspace}

Substantive part of the algorithm is drawing random ray based on available points and gradient vectors. We considered two variants of the algorithms differing in method to draw a subspace: vanilla, in which subspace is fully random, and cone-restricted, in which there is a dominant direction such that angle between random subspace and that direction is bounded. In both variants, $r$ is constructed as follows:
\begin{align*}
    r: b \to x_1 + A (b - x_1[B]),
\end{align*}
where $x[B]$ means vector from $\mathbb{R}^b$ composed of values of $x$ by indices $B$, and $A \in \mathbb{R}^{n \times b}$ is randomised matrix containing linear factors defining the affine subspace, which is different for variants under consideration.

In vanilla variant, $A$ is generated as follows:
\begin{align*}
    &c_{i j} \propto \mathcal{U}(-1, 1)\\
    &a_{i j} = \tan{(\pi/2 \cdot C)}\\
    &A[B] = I_b,
\end{align*}
where $I_b$ is $b \times b$ identity matrix, and $A[B]$ means submatrix of $A$ containing rows by indices $B$, analogously.

In cone-restricted variant, procedure is more complicated. Independently on dominant direction, angles under $\tan$ must be in $(-\pi/2, \pi/2)$. $A$, therefore, is generated as follows:
\begin{align*}
    &c_{i j} \propto \mathcal{U}(-1, 1)\\
    &\alpha_{i j} = \arctan{\frac{g_i}{g_j}}\\
    &\underline{\alpha}_{i j} = \max\{-\pi/2 + \epsilon, \alpha_{i j} - \beta \cdot a\}\\
    &\overline{\alpha}_{i j} = \min\{\alpha_{i j} + \beta \cdot a, \pi/2 - \epsilon\}\\
    &A = \tan{\left(\frac{\overline{\alpha} + \underline{\alpha}}{2} + \frac{\overline{\alpha} - \underline{\alpha}}{2} \circ C\right)}\\
    &A[B] = I_b,
\end{align*}
where $\epsilon = 10^{-16}$ is machine zero, $\beta$ is coefficient dependent on iteration of algorithm (can be used to increase angle), $a$ is initial angle of cone from which subspaces are drawn. $\circ$ denotes Hadamard product. $g \in \mathbb{R}^n$ is vector of dominant direction, it may be equal to $\nabla f(x_1)$ (first-order variant) or $x_2 - x_1$ (secant variant). 

\section{Numerical experiments} \label{sec:numerical}

To assess the practical performance of proposed algorithm and compare it with existing efficient algorithms, numerical experiments on quadratic (representing class of convex problems), Rosenbrock--Skokov (typical example of non-convex unimodal problem) and Rastrigin and DeVilliersGlasser02 functions (both are from benchmark \cite{gavana2018global}) are carried out.

Common details of further experiments are following. To solve the restricted problems of low dimensionality, Solar method use Nelder--Mead method \cite{nelder1965simplex} set to 10 iterations in all the cases. Each curve or point on presented plots is equipped with shadow: for classic randomised algorithms, sizes of lower and upper shadow are determined by standard deviation of function value measured in 3 or 5 runs with different random seeds; for Solar method, upper shadow is given by function value in run with worst final function value, and lower bound~-- with best final function value, correspondingly.

The implementation of algorithms used in the experiments and reference implementation of Solar method are in Python 3 and available as open source\footnote{\url{https://github.com/dmivilensky/Solar-method-non-convex-optimisation}}.

\subsection{Convex problem: quadratic function} \label{sec:convex}

Let's consider the problem of quadratic function optimisation with uniformly random positive-definite matrix $AA^\top$, some vector $b$ and scalar $c$:
\begin{align*}
    &f(x) = x^\top A A^\top x + b^\top x + c\\
    &a_{ij} \propto \mathcal{U}(0, a),\,b_i \propto \mathcal{U}(0, 1),\,c_i \propto \mathcal{U}(0, 1)
\end{align*}
with three particular instances of different dimensionality and conditioning:
\begin{center}
\begin{minipage}{3in}
\begin{enumerate}
    \item $a = 1,\,\kappa \approx 3.78 \cdot 10^3$\\
    $x \in \mathbb{R}^{10}, -20 \leq x \leq 20$\\
    $\|x^0 - x^*\|_2 \approx 23.16,\,f(x^*) \approx -4.77$
    \item $a = 5\sqrt{2},\,\kappa \approx 5.44 \cdot 10^4$\\
    $x \in \mathbb{R}^{25}, -5 \leq x \leq 5$\\
    $\|x^0 - x^*\|_2 \approx 8.40,\,f(x^*) \approx -1.42$
    \item $a = \sqrt{10},\,\kappa \approx 3.70 \cdot 10^5$\\
    $x \in \mathbb{R}^{50}, -5 \leq x \leq 5$\\
    $\|x^0 - x^*\|_2 \approx 1.23,\,f(x^*) \approx 0.71$
\end{enumerate}
\end{minipage}
\end{center}
where $\kappa = \lambda_{\max}(AA^\top) / \lambda_{\min}(AA^\top)$ is a condition number. In all the cases $\kappa \gg 1$, which makes chosen quadratic functions hard to optimise using standard gradient descent.

For this example, we apply Solar method in zeroth-order variant and compare it with zeroth-order methods of convex optimisation for fairness. We consider zeroth-order gradient descent with two-point feedback \cite{nesterov2017random} and its variant with line-search. To represent zeroth-order accelerated algorithms, we take Momentum three point method \cite{gorbunov2019stochastic} in deterministic setting.

\begin{figure}[H]
     \centering
     \begin{subfigure}[b]{0.32\textwidth}
         \centering
         \includegraphics[width=\textwidth]{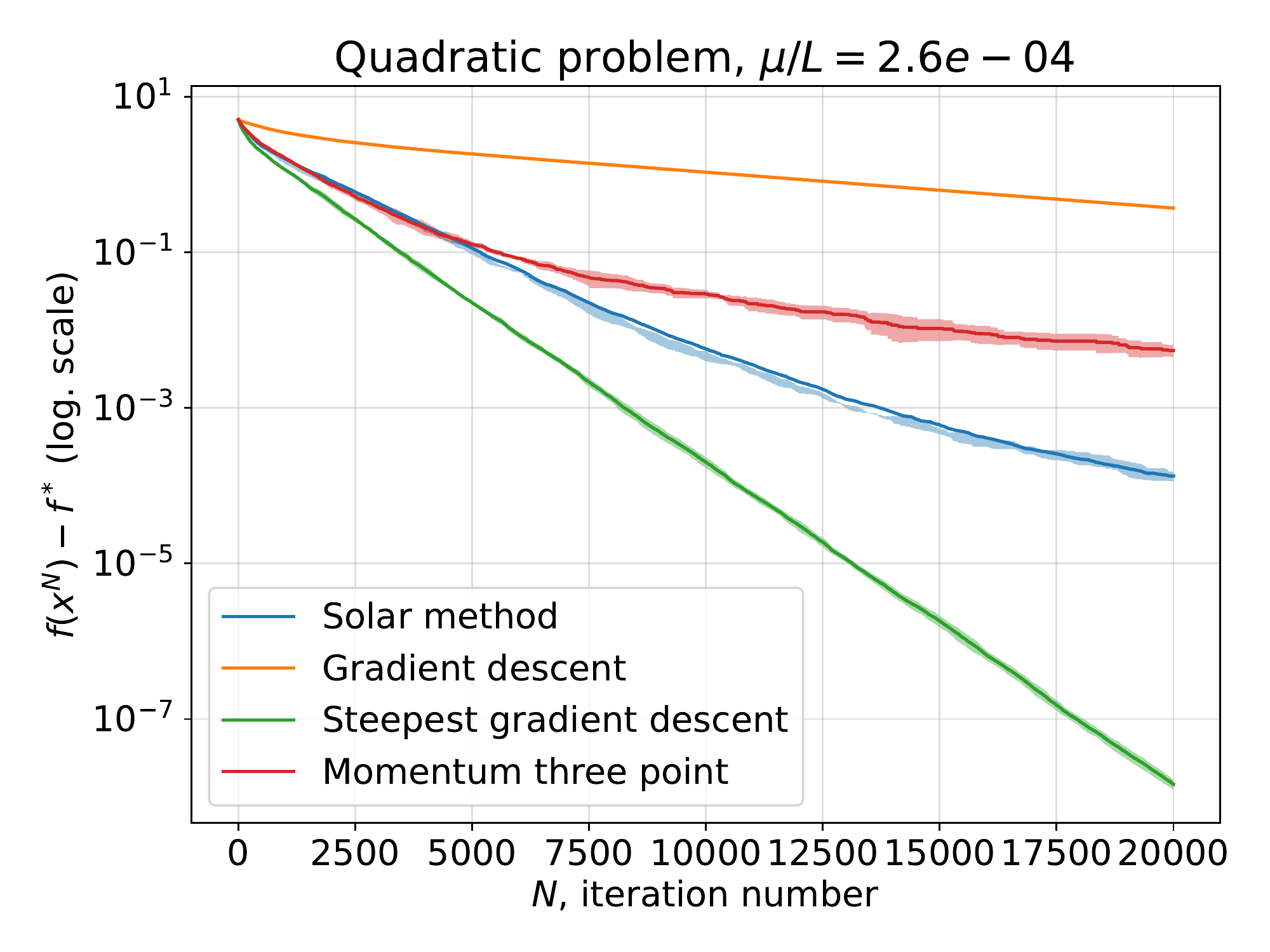}
         \caption{Convergence curves Solar method, two-point zeroth-order gradient descent with or without line-search and Momentum three point method}
        \label{fig:quadratic_a}
     \end{subfigure}
     \hspace{0.1cm}
     \begin{subfigure}[b]{0.32\textwidth}
         \centering
         \includegraphics[width=\textwidth]{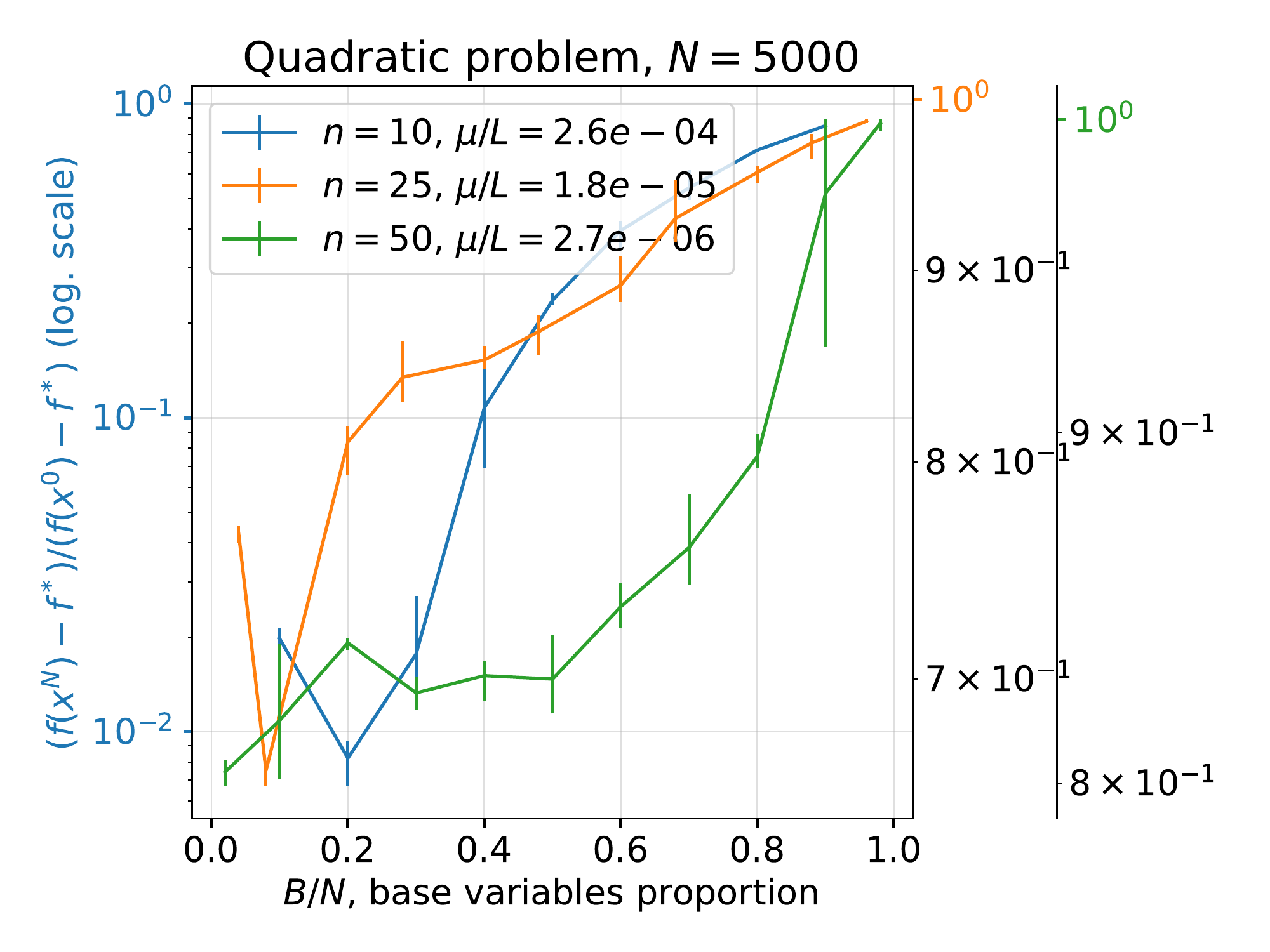}
         \caption{Function suboptimality achieved by Solar method on different quadratic functions with different proportion of base variables $B/N$}
        \label{fig:quadratic_b}
     \end{subfigure}
     \hspace{0.1cm}
     \begin{subfigure}[b]{0.32\textwidth}
         \centering
         \includegraphics[width=\textwidth]{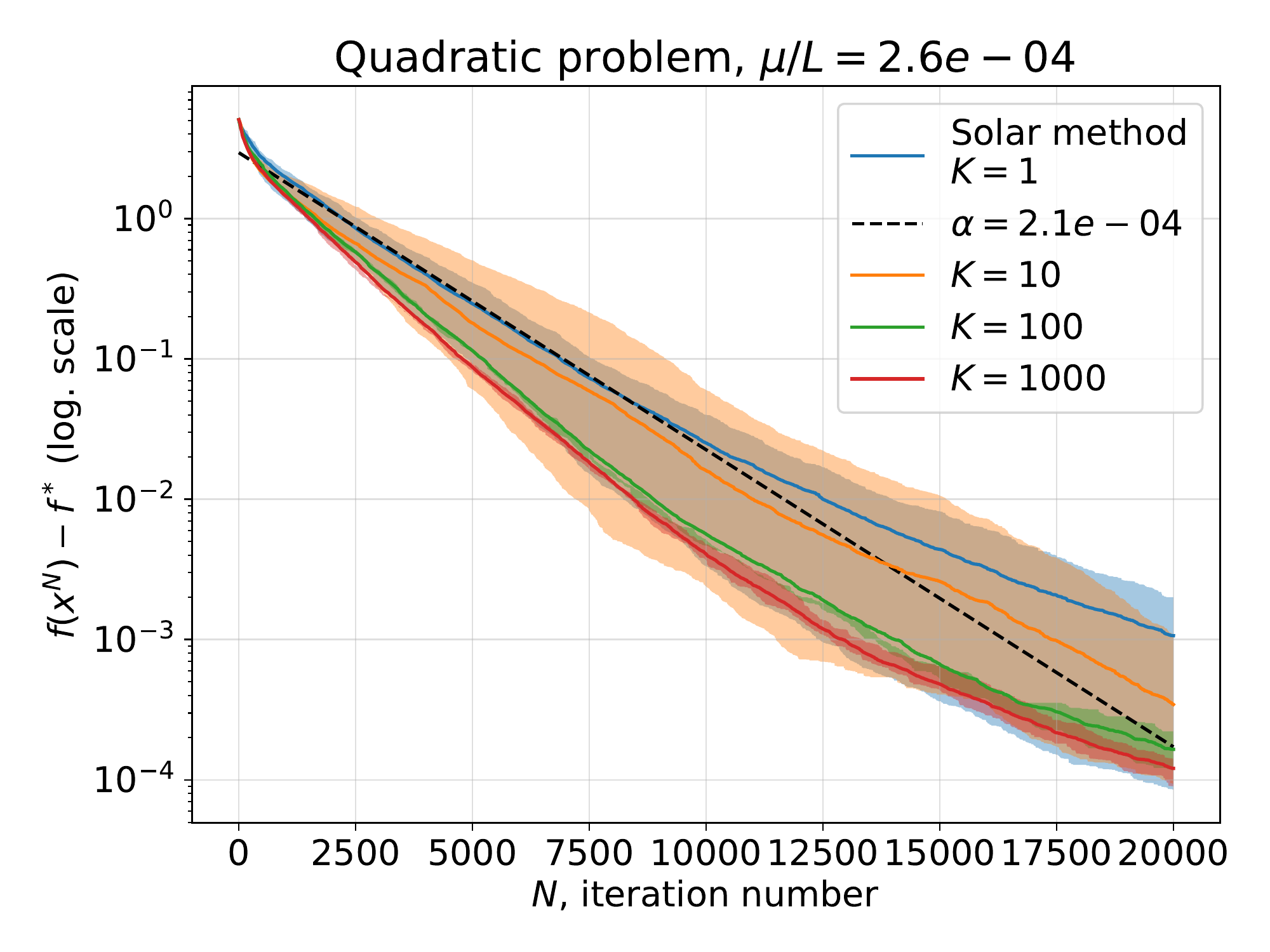}
         \caption{Convergence curves of Solar method for different $K=1, 10, 100, 1000$ and estimation of order of linear convergence}
        \label{fig:quadratic_c}
     \end{subfigure}
        \caption{Practical efficiency of Solar method on quadratic functions}
        \label{fig:quadratic}
\end{figure}

The description of the results presented on Figure~\ref{fig:quadratic} follows.

\textbf{Comparison of algorithms}. As is easily seen from Figure~\ref{fig:quadratic_a}, Solar method takes second place after steepest gradient descent. This is expectable that Solar method is not leading method for convex optimisation, because it puts much effort into exploration by choosing random direction instead of exploitation of given knowledge that gradient direction is descent direction leading to optimum. Unexpectedly, it performs better than Momentum three point method with theoretical convergence rate corresponding to accelerated methods.

\textbf{Dependency of convergence rate of problem's properties}. It is important to choose proper dimensionality $B$ of random subspaces, but it is expectable that optimal choice depends on properties of the problem, main of which are total dimensionality $N$ and conditioning $\mu/L \equiv 1/\kappa$. On Figure~\ref{fig:quadratic_b} convergence rate is estimated by relative function suboptimality after 5000 inner iterations of Solar method. One point represents convergence rate of method for given dimensionality, condition number and variable proportion of base variables $B/N$. It can be seen that form of dependence of convergence rate on dimensionality of subspaces varies, but as a rule low dimensionality is better. It can be explained by the uniform choice of number of iterations for method solving restricted problems, while the complexity of restricted problems grow with their dimensionality and more iterations are required. In some sense, curves reflect the dependency of error accumulated due to insufficient iterations of auxiliary method on dimensionality of subspaces. However, guided by limitation of budget, this plot should guide in choice of $B$. 

\textbf{Dependency of convergence rate on algorithm's hyperparameters}. Parameter for number of outer iterations $K$ of Solar method, or derived parameter for number of inner iterations per chosen set of base variables $N/K$, where $N$ is total number of iterations, does not seem to be essential. If random manifold is used for restricting the problem, this parameter can affect on variability of restricted problems, but linear subspaces are invariant on choice of base variables. Nevertheless, due to artefacts of random sampling (even in spherically symmetric sampling which we used) parameter $K$ affects practical performance. Figure~\ref{fig:quadratic_c} shows that the lower is number of inner iterations per chosen set of base variables, the less is dispersion of trajectories and the better resulting function value is on average. 

\textbf{Linear convergence in the beginning}. Another phenomenon Figure~\ref{fig:quadratic_c} shows is that in the beginning the convergence rate of Solar method is almost linear. We estimate order of convergence $\alpha$ by linear fit to show that it is close to $1/\kappa$, as for classic gradient descent.

\subsection{Unimodal problem: Rosenbrock--Skokov function} \label{sec:unimodal}

Let's consider multivariate generalisation of widely-known Rosenbrock non-convex test problem:
\begin{align*}
    &f(x) = (1 - x_1)^2 + 100 \cdot \sum_{i=2}^{100} \left(x_i - x_{i-1}^2\right)^2\\
    &x \in \mathbb{R}^{100}, -3 \leq x \leq 3\\
    &x^0 = (0.1,\,\dots,\,0.1)^\top, x^* = (1,\,\dots,\,1)^\top, f(x^*) = 0
\end{align*}

For this example, we consider variant of Solar method which uses gradient, together zeroth-order variant and compare them with first-order algorithms of conjugate gradient type. In particular, we run Fletcher--Reeves \cite{fletcher1964function} and Polak--Ribiere--Polyak \cite{polak1969note,polyak1969conjugate} conjugate gradient methods with or without restarts \cite{powell1977restart}. Besides, we test several variants of the Solar method differing in a way to draw random subspace on this problem.

\begin{figure}[H]
     \centering
     \begin{subfigure}[b]{0.33\textwidth}
         \centering
         \includegraphics[width=\textwidth]{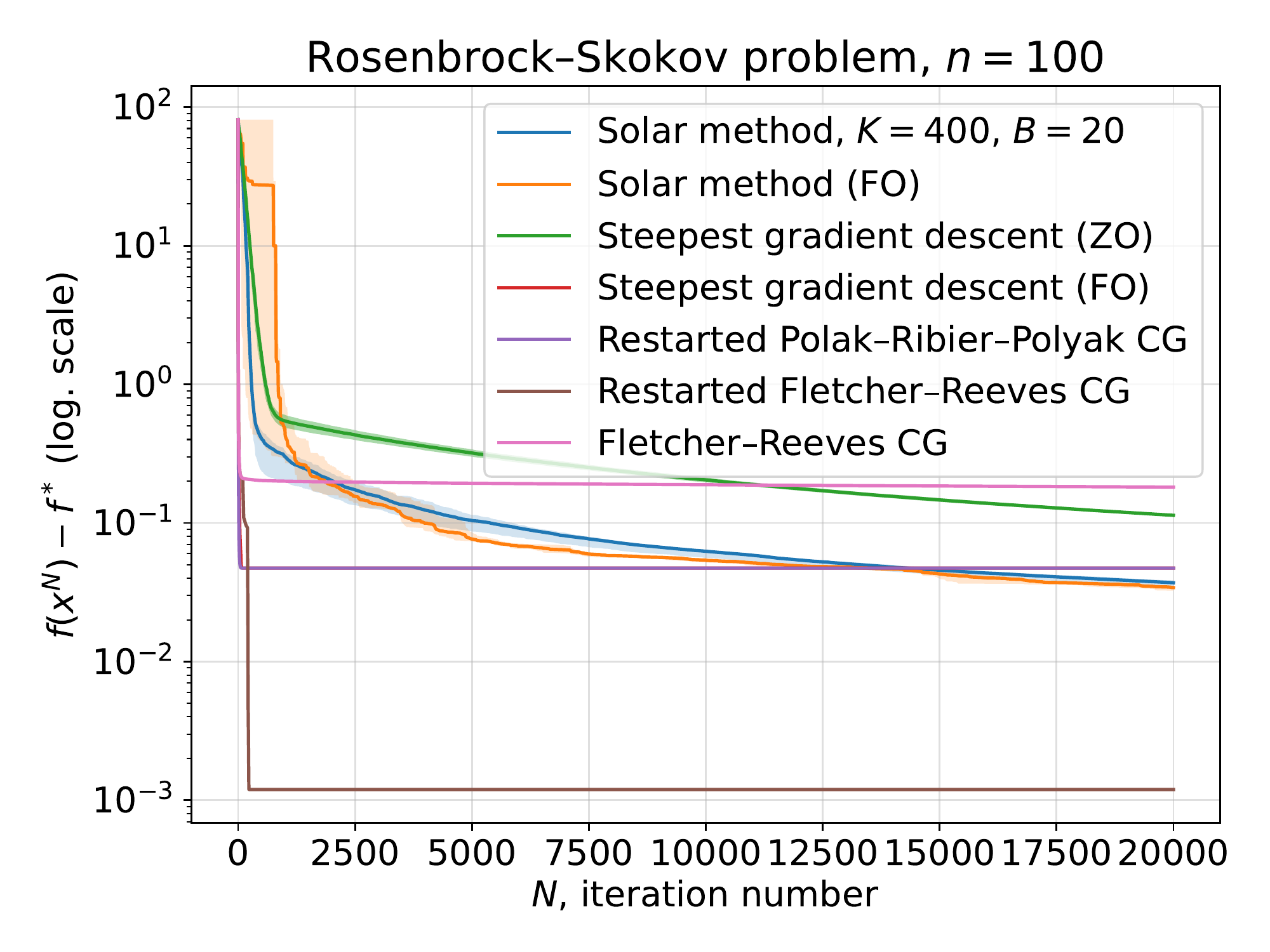}
         \caption{Convergence curves of Solar method, gradient descents, Fletcher--Reeves and Polak--Ribiere--Polyak conjugate gradients with or without restarts}
        \label{fig:rosenbrock_a}
     \end{subfigure}
     \hspace{0.1cm}
     \begin{subfigure}[b]{0.33\textwidth}
         \centering
         \includegraphics[width=\textwidth]{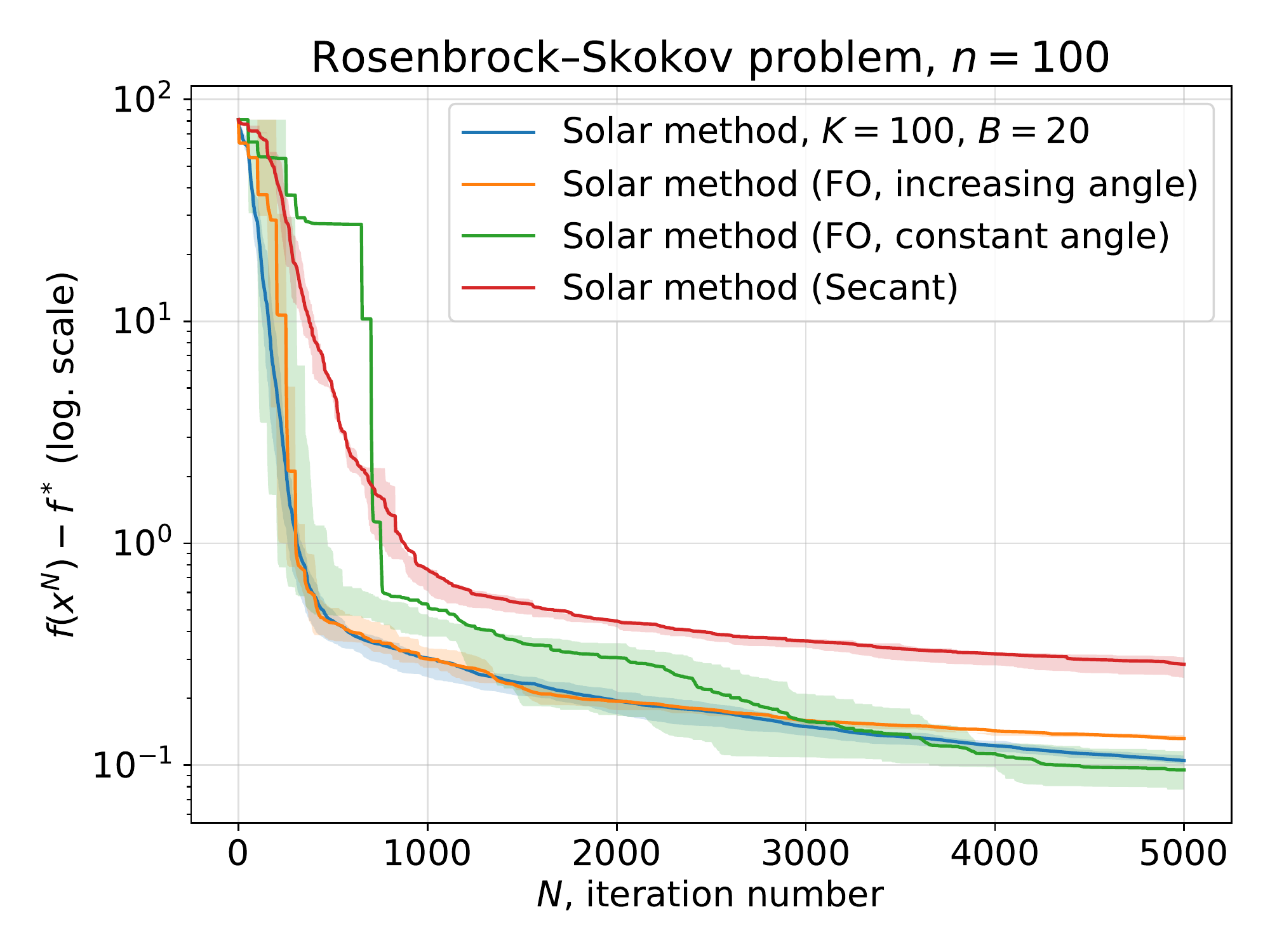}
         \caption{Convergence curves of several variants of Solar method with different ways to choose random subspace: sampling in cone around gradient or secant direction through two points}
        \label{fig:rosenbrock_b}
     \end{subfigure}
        \caption{Practical efficiency of Solar method on Rosenbrock--Skokov function}
        \label{fig:rosenbrock}
\end{figure}

The description of the results presented on Figure~\ref{fig:rosenbrock} follows.

\textbf{Comparison of algorithms}. As one can see from the Figure~\ref{fig:rosenbrock_a}, both variants of Solar method overtake Polak--Ribiere--Polyak algorithm and unrestarted Fletcher--Reeves algorithm (restarted Fletcher--Reeves algorithm is the only algorithm overtaking Solar method, and far ahead all other methods). Despite the instant convergence of conjugate gradient methods in the beginning, they stuck in a valley, while Solar method preserves steady convergence rate.

\textbf{Comparison of variants of Solar method}. There were tested four variants of Solar method: vanilla, in which subspace is defined by random linear dependency between base and the rest of variables, first-order, in which linear dependencies are in cone with constant angle around gradient, similar but with angle increasing over time and secant, in which two best points are maintained and cone is taken around segment connecting them. Briefly, the fastest variant is first-order one with constant angle, but the performance of all the variants is only slightly different, so in practice one should choose variant according to its efficiency on particular problem.

\subsection{Global problems: Rastrigin and DeVilliersGlasser02 functions} \label{sec:global}

Let's consider two global optimisation problems from benchmark \cite{gavana2018global} to assess the performance of Solar method in multi-extrema setting. First problem is for classical Rastrigin function, with 200 variables:
\begin{align*}
    &f(x) = 10\cdot200\cdot\sum_{i=1}^{200} (x_i^2 - 10 \cos{(2 \pi x_i)})\\
    &x \in \mathbb{R}^{200}, -5.12 \leq x \leq 5.12\\
    &x^0 = (5,\,\dots,\,5)^\top,\,x^* = (0,\,\dots,\,0)^\top,\,f(x^*) = 0
\end{align*}
and second is the hardest problem in \cite{gavana2018global}, for DeVilliersGlasser02 function of only 5 variables:
\begin{align*}
    &f(x) = \sum_{i=1}^{24} \left(x_1 x_2^{t_i} \tanh{(x_3 t_i + \sin{(x_4 t_i)})} \cos{(t_i e^{x_5})} - y_i\right)^2,\\
    &\text{where }t_i = 0.1 (i - 1),\,y_i = 53.81 \cdot 1.27^{t_i} \tanh{(3.012 t_i + \sin{(2.13 t_i)})} \cos{(t_i e^{0.507})}\\
    &x \in \mathbb{R}^{5}, 1 \leq x \leq 60\\
    &x^0 = (30,\,\dots,\,30)^\top,x^* = (53.81, 1.27, 3.012, 2.13, 0.507)^\top,\,f(x^*) = 0
\end{align*}

For this example, we compare Solar method with two classic and practically efficient algorithms of global optimisation: Simulated Annealing \cite{xue1993parallel}, and Monotonic Sequence Basin Hopping \cite{leary2000global}.

\begin{figure}[H]
     \centering
     \begin{subfigure}[b]{0.33\textwidth}
         \centering
         \includegraphics[width=\textwidth]{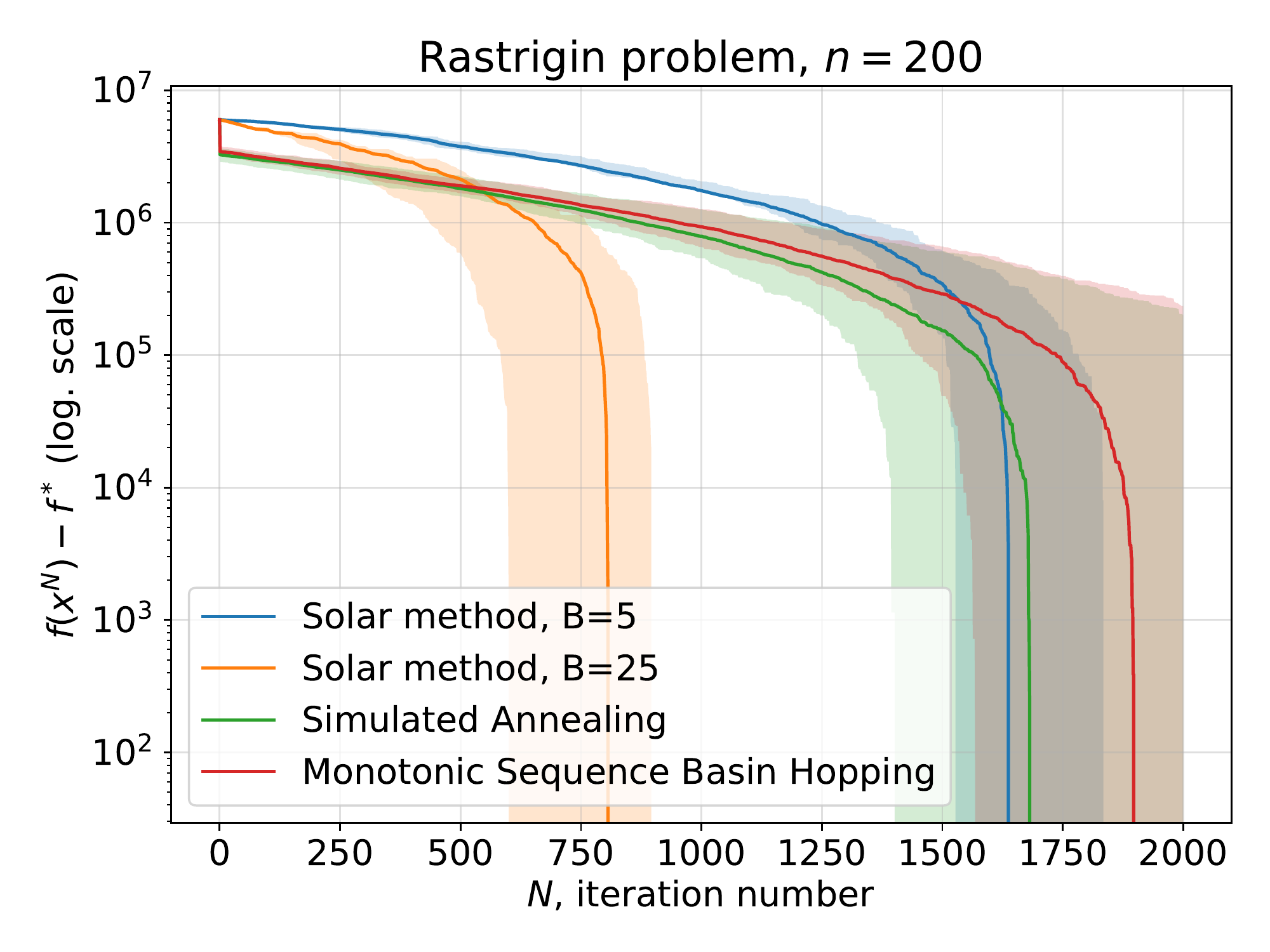}
         \caption{Convergence curves of Solar method, Simulated Annealing and MSBH on Rastrigin problem}
     \end{subfigure}
     \hspace{0.1cm}
     \begin{subfigure}[b]{0.33\textwidth}
         \centering
         \includegraphics[width=\textwidth]{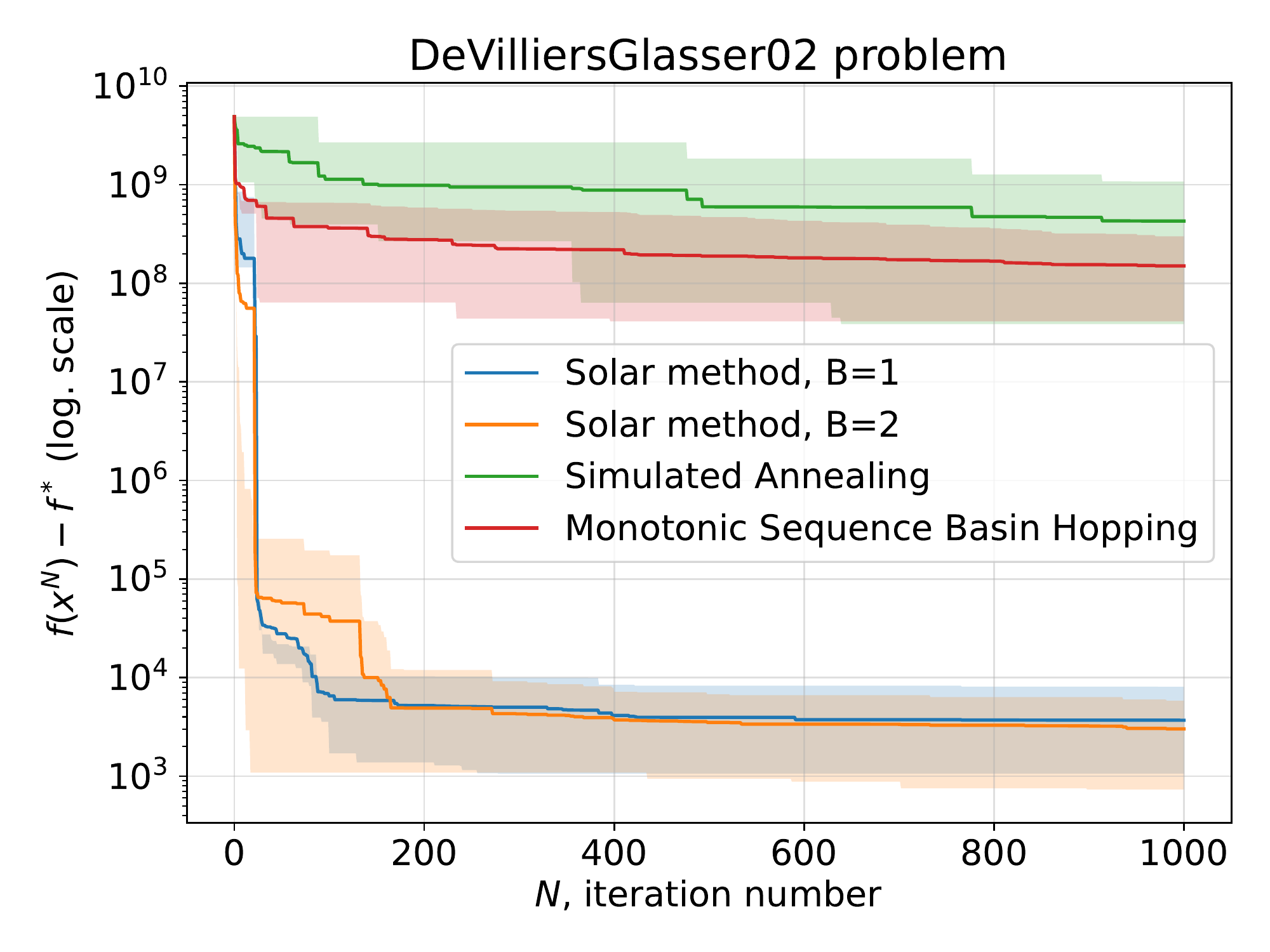}
         \caption{Convergence curves of Solar method, Simulated Annealing and MSBH on DeVilliersGlasser02 problem}
     \end{subfigure}
        \caption{Practical efficiency of Solar method on global optimisation test functions}
        \label{fig:global}
\end{figure}

The results presented on Figure~\ref{fig:global} can be described together as both of experiments compare algorithms' performance. In both cases Solar method turns out to be more efficient in search of deeper local optimum. Note, that for better performance dimensionality of subspaces should be chosen not very small. The explanation of the computational fact that Solar method is more efficient than standard hopping algorithms is that Solar method is less sensitive to choice of hyperparameters. There is no doubt that after grid search of jump length and proper temperature decreasing policy classic algorithms will achieve at least not worse function value in the same time. But this search of hyperparameters is time-consuming. In turn, Solar method has a great advantage that allows it jump farther: instead of localisation of the jump by limiting its length, Solar method restricts it on low-dimensional subspace, which does not take away the opportunity to find a better local optimum far from current best point.

\section{Discussion}

This paper proposes heuristic global optimisation algorithm called Solar method, which is based on optimisation over randomly drawn subspaces to search for a point step to. The algorithm demonstrated competitive or leading convergence rate on convex, unimodal and global test functions in comparison with algorithms known as the best in practice general-purpose methods in corresponding classes of problems. 

The proposed algorithm is an instance of more general scheme, so its further development shall be aimed at more delicate exploitation of advantages of the scheme itself. Theoretical and empirical study of this scheme is one barely permeable but, we believe, very perspective direction for the future work.

From the practical point of view, the usability of the algorithm itself shall be improved. Firstly, more detailed study of how the inexactness in solving the auxiliary problems accumulate and affect on overall performance is required. Secondly, current technical solution of maintaining the $k$ best points is too time-consuming, so one needs to find more convenient data structure for this purpose.

Finally, algorithm can be extended in several directions. Firstly, it is needed to explore more options of choice of random subspace and consider random manifold generality. Secondly, the algorithm can use itself as an algorithm to solve restricted problems, which gives complex multi-level algorithm with unknown performance and properties. 

\bibliographystyle{plain}
\bibliography{main}

\begin{thebibliography}{10}

\bibitem{cartis2022randomised}
Coralia Cartis, Jaroslav Fowkes, and Zhen Shao.
\newblock Randomised subspace methods for non-convex optimization, with
  applications to nonlinear least-squares.
\newblock {\em arXiv preprint arXiv:2211.09873}, 2022.

\bibitem{fletcher1964function}
Reeves Fletcher and Colin~M. Reeves.
\newblock Function minimization by conjugate gradients.
\newblock {\em The computer journal}, 7(2):149--154, 1964.

\bibitem{gavana2018global}
Andrea Gavana.
\newblock Global optimization benchmarks and ampgo, 2018.

\bibitem{gorbunov2019stochastic}
Eduard Gorbunov, Adel Bibi, Ozan Sener, El~Houcine Bergou, and Peter
  Richt{\'a}rik.
\newblock A stochastic derivative free optimization method with momentum.
\newblock {\em arXiv preprint arXiv:1905.13278}, 2019.

\bibitem{gower2019rsn}
Robert Gower, Dmitry Kovalev, Felix Lieder, and Peter Richt{\'a}rik.
\newblock R{S}{N}: Randomized subspace {N}ewton.
\newblock {\em Advances in Neural Information Processing Systems}, 32, 2019.

\bibitem{gower2021stochastic}
Robert~M. Gower, Peter Richt{\'a}rik, and Francis Bach.
\newblock Stochastic quasi-gradient methods: Variance reduction via {J}acobian
  sketching.
\newblock {\em Mathematical Programming}, 188:135--192, 2021.

\bibitem{grishchenko2021proximal}
Dmitry Grishchenko, Franck Iutzeler, and J{\'e}r{\^o}me Malick.
\newblock Proximal gradient methods with adaptive subspace sampling.
\newblock {\em Mathematics of Operations Research}, 46(4):1303--1323, 2021.

\bibitem{jones2001taxonomy}
Donald~R. Jones.
\newblock A taxonomy of global optimization methods based on response surfaces.
\newblock {\em Journal of global optimization}, 21:345--383, 2001.

\bibitem{kozak2021stochastic}
David Kozak, Stephen Becker, Alireza Doostan, and Luis Tenorio.
\newblock A stochastic subspace approach to gradient-free optimization in high
  dimensions.
\newblock {\em Computational Optimization and Applications}, 79(2):339--368,
  2021.

\bibitem{leary2000global}
Robert~H. Leary.
\newblock Global optimization on funneling landscapes.
\newblock {\em Journal of Global Optimization}, 18(4):367, 2000.

\bibitem{nelder1965simplex}
John~A. Nelder and Roger Mead.
\newblock A simplex method for function minimization.
\newblock {\em The computer journal}, 7(4):308--313, 1965.

\bibitem{nemirovskij1983problem}
Arkadij~S. Nemirovskij and David~B. Yudin.
\newblock {\em Problem complexity and method efficiency in optimization}.
\newblock Wiley-Interscience, 1983.

\bibitem{nesterov2017random}
Yurii Nesterov and Vladimir Spokoiny.
\newblock Random gradient-free minimization of convex functions.
\newblock {\em Foundations of Computational Mathematics}, 17:527--566, 2017.

\bibitem{nesterov2012efficiency}
Yurii~E. Nesterov.
\newblock Efficiency of coordinate descent methods on huge-scale optimization
  problems.
\newblock {\em SIAM Journal on Optimization}, 22(2):341--362, 2012.

\bibitem{pasechnyuksolar}
Dmitry~A. Pasechnyuk and Alexander Gornov.
\newblock Solar method for non-convex problems: hypothesizing approach to an
  optimization.
\newblock Retrieved from
  \url{http://dmivilensky.ru/preprints/Solar\%20method\%20for\%20non-convex\%20problems.pdf},
  2022.

\bibitem{patrascu2015efficient}
Andrei Patrascu and Ion Necoara.
\newblock Efficient random coordinate descent algorithms for large-scale
  structured nonconvex optimization.
\newblock {\em Journal of Global Optimization}, 61(1):19--46, 2015.

\bibitem{polak1969note}
Elijah Polak and Gerard Ribiere.
\newblock Note sur la convergence de m{\'e}thodes de directions conjugu{\'e}es.
\newblock {\em Revue fran{\c{c}}aise d'informatique et de recherche
  op{\'e}rationnelle. S{\'e}rie rouge}, 3(16):35--43, 1969.

\bibitem{polyak1969conjugate}
Boris~T. Polyak.
\newblock The conjugate gradient method in extremal problems.
\newblock {\em USSR Computational Mathematics and Mathematical Physics},
  9(4):94--112, 1969.

\bibitem{powell1977restart}
Michael James~David Powell.
\newblock Restart procedures for the conjugate gradient method.
\newblock {\em Mathematical programming}, 12:241--254, 1977.

\bibitem{richtarik2014iteration}
Peter Richt{\'a}rik and Martin Tak{\'a}{\v{c}}.
\newblock Iteration complexity of randomized block-coordinate descent methods
  for minimizing a composite function.
\newblock {\em Mathematical Programming}, 144(1-2):1--38, 2014.

\bibitem{torn1989global}
Aimo T{\"o}rn and Antanas Zilinskas.
\newblock {\em Global optimization}, volume 350.
\newblock Springer, 1989.

\bibitem{wales1997global}
David~J. Wales and Jonathan~P.K. Doye.
\newblock Global optimization by basin-hopping and the lowest energy structures
  of lennard--jones clusters containing up to 110 atoms.
\newblock {\em The Journal of Physical Chemistry A}, 101(28):5111--5116, 1997.

\bibitem{xu2017globally}
Yangyang Xu and Wotao Yin.
\newblock A globally convergent algorithm for nonconvex optimization based on
  block coordinate update.
\newblock {\em Journal of Scientific Computing}, 72(2):700--734, 2017.

\bibitem{xue1993parallel}
Guo-Liang Xue.
\newblock Parallel two-level simulated annealing.
\newblock In {\em Proceedings of the 7th international conference on
  Supercomputing}, pages 357--366, 1993.

\bibitem{zhigljavsky2012theory}
Anatoly~A. Zhigljavsky.
\newblock {\em Theory of global random search}, volume~65.
\newblock Springer Science \& Business Media, 2012.

\end{thebibliography}

\end{document}